# Inverse applications of the generalized Littlewood theorem concerning integrals of the logarithm of analytic functions: an easy method to establish equalities between different analytic functions

**Sergey K. Sekatskii**[1,*]

[1] Laboratory of Biological Electron Microscopy, IPHYS, BSP 419, Ecole Polytechnique Fédérale de Lausanne, and Department of Fundamental Biology, Faculty of Biology and Medicine, University of Lausanne, CH1015 Lausanne, Switzerland.

[*] Correspondence: Serguei.Sekatski@epfl.ch; Tel.: +41-21-693-0445

**ABSTRACT**

Recently, we have established and used the generalized Littlewood theorem concerning contour integrals of the logarithm of analytical function to obtain new criteria equivalent to the Riemann hypothesis. Later, the same theorem was applied to calculate certain infinite sums and study the properties of zeroes of a few analytical functions. In this Note, we discuss what, in a sense, are inverse applications of this theorem. We first prove an easy Lemma that if two meromorphic on the whole complex plane functions $f(z)$ and $g(z)$ have the same zeroes and poles, taking into account their orders, and have appropriate asymptotic for large $|z|$, then for some integer $n$ $\frac{d^n \ln(f(z))}{dz^n} = \frac{d^n \ln(g(z))}{dz^n}$. The use of this Lemma enables easy proofs of many identities between elliptic functions, their transformation and $n$-tuple rules. In particular, we show how for any complex number $a$, $\wp(z)-a$, where $\wp(z)$ is Weierstrass $\wp$-function, can be presented as a product and ratio of three elliptic $\theta_1$-functions of certain arguments. We also establish $n$-tuple rules for elliptic theta-functions and $\wp_z(z)$ derivatives.

# 1. Introduction

The generalized Littlewood theorem concerning contour integrals of the logarithm of analytical function is stated as follows [1, 2]:

**Theorem 1 (The Generalized Littlewood theorem):** *Let C denote the rectangle bounded by the lines $x = X_1$, $x = X_2$, $y = Y_1$, $y = Y_2$ where $X_1 < X_2$, $Y_1 < Y_2$ and let f(z) be analytic and non-zero on C and meromorphic inside it, and let also g(z) be analytic on C and meromorphic inside it. Let F(z)=ln(f(z)) be the logarithm defined as follows: we start with a particular determination on $x = X_2$, and obtain the value at other points by continuous variation along y=const from $\ln(X_2 + iy)$. If, however, this path would cross a zero or pole of f(z), we take F(z) to be $F(z \pm i0)$ according as to whether we approach the path from above or below. Let also $\tilde{F}(z) = \ln(f(z))$ be the logarithm defined by continuous variation along any smooth curve fully lying inside the contour which avoids all poles and zeroes of f(z) and starts from the same particular determination on $x = X_2$. Suppose also that the poles and zeroes of the functions f(z), g(z) do not coincide.*

*Then*

$$\int_C F(z)g(z)dz = 2\pi i \left( \sum_{\rho_g} res(g(\rho_g) \cdot \tilde{F}(\rho_g)) - \sum_{\rho_f^0} \int_{X_1 + iY_\rho^0}^{X_\rho^0 + iY_\rho^0} g(z)dz + \sum_{\rho_f^{pol}} \int_{X_1 + iY_\rho^{pole}}^{X_\rho^{pole} + iY_\rho^{pole}} g(z)dz \right) \quad (1)$$

*where the sum is over all $\rho_g$ which are poles of the function g(z) lying inside C, all $\rho_f^0 = X_\rho^0 + iY_\rho^0$ which are zeroes of the function f(z) both counted taking into account their multiplicities (that is the corresponding term is multiplied by m for a zero of the order m) and which lie inside C, and all $\rho_f^{pole} = X_\rho^{pole} + iY_\rho^{pole}$ which are poles of the function f(z) counted taking into account their multiplicities and which lie inside C. The assumption is that all relevant integrals on the right hand side of the equality exist.*

The proof of this theorem [2] is very close to the proof of the standard Littlewood theorem corresponding to the case *g(z)=1*, see e.g. [3]. Especially interesting are some particular cases when the contour integral $\int_C F(z)g(z)dz$ disappears (tends to zero) if the contour tends to infinity, that is when $X_1, Y_1 \to -\infty, X_2, Y_2 \to +\infty$. (More precisely, when we can find some sequence of the contours $C_j$ tending to infinity and such, that the contour integral tends to zero, see below). This means that eq. (1) takes the form

$$\sum_{\rho_f^0} \int_{-\infty+iY_\rho^0}^{X_\rho^0+iY_\rho^0} g(z)dz - \sum_{\rho_f^{pol}} \int_{-\infty+iY_\rho^{pole}}^{X_\rho^{pole}+iY_\rho^{pole}} g(z)dz = \sum_{\rho_g} res(g(\rho_g) \cdot F(\rho_g)) \qquad (2)$$

If the integrals here can be calculated explicitly, in this way one obtains equalities involving finite or infinite sums (this last case is the most interesting one). Of course, indeed this is not necessary that the integrals $\int_C F(z)g(z)dz$ disappear for *an arbitrary* sequence of contours $C_i$ tending to infinity; this is enough to find *some specific* sequence such that the integrals do disappear in the limit. Note also, that the orthogonal parallelograms described in the conditions of the Theorem 1 actually are not necessary and much broader family of the contours can be considered.

Earlier, this approach was used by us to analyze some properties of the zeroes of the Riemann zeta-function (see e.g. [4] for a general discussion of this function) - in particular, to establish a number of theorems equivalent to the Riemann hypothesis [1, 2, 5, 6]. (Some of these results were recently included in the corresponding chapter of the Encyclopedia of Mathematics and its Applications [7]). In Ref. [8], we discuss the use of the Generalized Littlewood Theorem to calculate many infinite sums over integers and to study the properties of zeroes of some analytical functions, such as incomplete Riemann and gamma-functions and polygamma functions. In Ref. [9] this approach with the same purposes has been applied to elliptic functions.

In the present paper we discuss what, in a sense, is the *inverse* application of this generalized Littlewood Theorem. Namely, we demonstrate how from the circumstance that certain analytical functions have the same poles and zeroes, together with some additional information, the formulae connecting these functions can be established. These applications are first illustrated by certain trigonometrical functions and gamma-function, and then applied to easily prove numerous equalities and transformation rules between different elliptic functions.

## 2. The main Lemma

Let us state the following easy Lemma 1.

**Lemma 1.** *Let the functions f(z) and g(z) be analytic and meromorphic on the complex plane, and let for some integer n there exists a sequence of contours $C_i$ tending to infinity, such as defined in the conditions of the generalized Littlewood theorem, and such that $\int_{C_i} \frac{1}{(z-a)^{n+1}} \ln(f(z))dz$ and $\int_{C_i} \frac{1}{(z-a)^{n+1}} \ln(g(z))dz$ tend to zero. Here a is an arbitrary complex number not coinciding with any zero or pole of the functions f(z) and g(z). Let also the poles and zeroes of these functions, taking into account their multiplicities, coincide. Then $\frac{d^n \ln(f(z))}{dz^n} = \frac{d^n \ln(g(z))}{dz^n}$.*

**Proof.** Actually, there is nothing to prove because we just restate the conditions of the generalized Littlewood theorem. The contour integrals $\int_{C_i} \frac{1}{(z-a)^{n+1}} \ln(f(z))dz$ and $\int_{C_i} \frac{1}{(z-a)^{n+1}} \ln(g(z))dz$ for contours $C_i$ and points *a* considered above tend to zero, and the identity of poles and zeroes of the functions f(z), g(z) takes place. Thus, according to eq. (2), we have $\frac{d^n \ln(f(z))}{dz^n}|_{z=a} = \frac{d^n \ln(g(z))}{dz^n}|_{z=a}$ for all $a$ not coinciding with any zero or pole of the functions $f(z)$, $g(z)$. □

Certainly, this Lemma has close connection with the famous Cauchy - Liouville's theorem stating that any bounded entire analytical function is a constant, and its generalizations by Hadamard and others, see e.g. [3]. Nevertheless, it turns out to be rather useful in applications and thus, in our opinion, deserves a special attention.

## 3. Illustrative applications
Let us first consider some simple illustrations.
### 3.1 Gamma-function
First, we apply Lemma 1 with $n=2$ to the functions $\Gamma(2z)$ and $\Gamma(z)\Gamma(z+1/2)$. Both these functions have the same poles at $z=0, -1/2, -1, -3/2, -2...$ and no zeroes, and for large $|z|$ the asymptotic of their logarithms is $O(z \ln z)$; see e.g. [10, 11] for the general discussion of gamma-function and its derivatives. Thus Lemma 1 gives $\psi'(z) + \psi'(z+1/2) = 4\psi'(2z)$. Integration gives $\psi(z) + \psi(z+1/2) = 2\psi(2z) + C_1$, and one more integration gives $\Gamma(z)\Gamma(z+1/2) = C_2 e^{C_1 z} \Gamma(2z)$, where $C_1$ and $C_2$ are integration constants to be determined.

Now the elementary $\Gamma(1) = \Gamma(2) = 1$, $\Gamma(1/2) = \sqrt{\pi}$ and $\Gamma(3/2) = \Gamma(1+1/2) = \frac{1}{2}\sqrt{\pi}$ enable, after the substitution of $z=1/2$ and $z=1$ into above equality, to recover the Legendre duplication rule:

$$\Gamma(z)\Gamma(z+1/2) = 2^{1-2z} \sqrt{\pi} \Gamma(2z) \qquad (3).$$

The reflection formula

$$\Gamma(1-z)\Gamma(z) = \frac{\pi}{\sin(\pi z)} \qquad (4)$$

is also trivially obtained by applying Lemma 1 with $n=2$ to $\Gamma(1-z)\Gamma(z)$ and $\frac{1}{\sin(\pi z)}$.

Thus $\psi'(1-z) + \psi'(z) = -\frac{\pi^2}{\cos(\pi z)}$. The first integration gives

$-\psi(1-z) + \psi(z) = -\pi\cot(\pi z) + C_1$; substitution of $z=1/2$, without even the knowledge of $\psi(1/2)$, gives $C_1 = 0$. Second integration gives $\Gamma(1-z)\Gamma(z) = \frac{C_2}{\sin(\pi z)}$, substitution of $z=1/2$ gives $C_2 = \pi$.

Let us now prove the same formula in slightly different way relaying on the Laurent expansions of the corresponding logarithms. We know $\Gamma^{-1}(1+z) = 1 + \gamma z + O(z^2)$ so that from the functional equation $\Gamma(1+z) = z\Gamma(z)$ we have $\Gamma^{-1}(z) = z + \gamma z^2 + O(z^3)$, $\Gamma^{-1}(z)\Gamma(1-z) = z + O(z^3)$, and $\ln(\Gamma^{-1}(z)\Gamma(1-z)/z) = O(z^2)$. We also know that $\ln(\frac{\sin \pi z}{\pi z}) = O(z^2)$. Generalized Littlewood theorem (its inverse application, Lemma 1) shows that $O(z^2)$ terms in both these formulae are identical. Thus $\ln(\Gamma^{-1}(z)\Gamma^{-1}(1-z)/z) = \ln\frac{\sin \pi z}{\pi z}$ and (4) follows.

In exactly the same way we obtain an easy proof of $\sum_{k=0}^{n-1} \psi'(z + \frac{k}{n}) = n^2 \psi'(nz)$ and, by integration, $\sum_{k=0}^{n-1} \psi(z + \frac{k}{n}) = n\psi(nz) + C_1(n)$. Substituting $z=1$ we get $\psi(1) + \sum_{k=1}^{n-1} \psi(1 + \frac{k}{n}) = n\psi(n) + C_1(n)$. The use of $\psi(1) = -\gamma$ and $\psi(1+z) = \psi(z) + \frac{1}{z}$, whence $\psi(n) = \sum_{k=1}^{n-1}\frac{1}{k} - \gamma$, gives $-\gamma + \sum_{k=1}^{n-1}\psi(\frac{k}{n}) + n\sum_{k=1}^{n-1}\frac{1}{k} = n(\sum_{k=1}^{n-1}\frac{1}{k} - \gamma) + C_1(n)$, so that $C_1(n) = \sum_{k=1}^{n-1} \psi(\frac{k}{n}) + (n-1)\gamma$. Applying Gauss' identity $\sum_{k=1}^{n} \psi(\frac{k}{n}) = -n(\gamma + \ln n)$ [10, 11], we obtain $C_1(n) = -n\ln n$ whence $\sum_{k=0}^{n-1} \psi(z + \frac{k}{n}) = n\psi(nz) - n\ln n$. One more integration gives $\prod_{k=0}^{n-1} \Gamma(z + \frac{k}{n}) = C_2(n) e^{-nz \ln n} \Gamma(nz)$ and $C_2(n) = n \prod_{k=1}^{n-1} \Gamma(\frac{k}{n})$, but unfortunately we do not see an obvious way to find that $C_2(n) = (2\pi)^{\frac{n-1}{2}} n^{\frac{1}{2}}$ for the case, and thus fully restore the Gauss multiplication theorem [10, 11]:

$$\prod_{k=0}^{n-1} \Gamma(z + \frac{k}{n}) = (2\pi)^{\frac{n-1}{2}} n^{\frac{1}{2}} e^{-nz \ln n} \Gamma(nz) \qquad (5).$$

*3.2. Trigonometrical functions*

The corresponding relations for trigonometrical functions are very well known, so we limit ourselves with quite short exposition (it also will be useful when similar transformation rules for elliptical functions will be discussed; see below). The equivalence of the sums over poles and zeroes expressed via $\int_C \frac{1}{(z-a)^3} \ln(\sin(nz)) dz$ and $\int_C \frac{1}{(z-a)^3} \ln(\prod_{k=0}^{n-1} \sin(z + \frac{\pi k}{n})) dz$ (Lemma 1) readily furnishes

$$\sum_{k=0}^{n-1} \frac{1}{\sin^2(z + \frac{\pi k}{n})} = \frac{n^2}{\sin^2(nz)} \qquad (6).$$

The integration of (6) gives $\sum_{k=0}^{n-1} \cot(z + \frac{\pi k}{n}) = n\cot(nz) + C(n)$, the substitution $z = \frac{\pi}{2n}$ shows that $C(n)=0$ whence

$$\sum_{k=0}^{n-1} \cot(z + \frac{\pi k}{n}) = \cot(nz) \qquad (7).$$

The next integration gives

$$\prod_{k=0}^{n-1} \sin(z + \frac{k\pi}{n}) = C_1(n) \sin(nz) \qquad (8).$$

This is a well-known (and very good) students' exercise to establish Euler relation $\prod_{k=1}^{n-1}\sin(\frac{k\pi}{n}) = \frac{n}{2^{n-1}}$. This relation corresponds to the case when $z$ tends to zero in eq. (8). Thus

$$\prod_{k=0}^{n-1}\sin(z + \frac{k\pi}{n}) = 2^{1-n}\sin(nz) \qquad (9).$$

Finally, let us note the following rarely mentioned equality between the products and sums of squares for the sine function. From (9) and (6) it immediately follows that

$$\prod_{k=0}^{n-1}\sin^2(z + \frac{k\pi}{n}) = n^2 2^{2-2n} \frac{1}{\sum_{k=0}^{n-1}\sin^{-2}(z + \frac{k\pi}{n})} \qquad (10).$$

Attention on the sign should be paid if extracting the square root of (10).

## 4. Applications to elliptic functions
*4.1 Definitions and main properties of elliptic functions*

Elliptic functions, which are much studied because of their high importance for mathematics and physics, see e.g. [12-18] (whenever possible, below we will mostly cite "encyclopedia-like" Ref. [12]), are a fertile ground for the approach based on the generalized Littlewood theorem. First we need to give necessary definitions and specify a notation used, especially because, unfortunately, different notations and conventions still co-exist in the elliptic functions research field. We define four theta-functions as follows

$$\theta_1(z,\ q) = \theta_1(z|\tau) = 2\sum_{n=0}^{\infty}(-1)^n q^{(n+1/2)^2}\sin((2n+1)z) \qquad (11a)$$

$$\theta_2(z,\ q) = \theta_2(z|\tau) = 2\sum_{n=0}^{\infty}q^{(n+1/2)^2}\cos((2n+1)z) \qquad (11b)$$

$$\theta_3(z,\ q) = \theta_3(z|\tau) = 1 + 2\sum_{n=0}^{\infty}q^{n^2}\cos(2nz) \qquad (11c)$$

$$\theta_4(z, q) = \theta_4(z|\tau) = 1 + 2\sum_{n=0}^{\infty}(-1)^n q^{n^2}\cos(2nz) \qquad (11d).$$

Here $q = e^{i\pi\tau}$ (it is named a nome), and $\text{Im}\,\tau > 0$. As functions of $z$ for any fixed $\tau$, they are entire and $2\pi$ - periodic, and they are quasiperiodic on the lattice formed by the points $z_{m,n} = (m + n\tau)\pi$ - in the sense, that the following relation holds:

$$\theta_1(z + (m + n\tau)\pi \mid \tau) = (-1)^{m+n} q^{-n^2} e^{-2inz}\theta_1(z|\tau) \qquad (12),$$

and similar relations exist for other theta-functions. Here and below, unless specifically stated to the contrary, $n, m \in Z$, and we will not always repeat this statement. The notation $\theta_j(z, q)$ and $\theta_j(z|\tau)$ is used on the equal footing.

These properties and defined by them asymptotic guaranty that for large $|z|$ the theta-function is at most $O(\exp(C|z|^2))$ with some constant $C$), hence in the limit of infinitely large contours we have the disappearance of the contour integrals $\int_{C_i}\frac{1}{(z-a)^{n+1}}\ln(\theta_j(z))dz$ for $j=1, 2, 3, 4$ and $n=3, 4\ldots$ The location of zeroes $\rho_i$, which all are simple, for these functions, also is well known. Namely, the functions $\theta_j(z|\tau)$ for $j=1, 2, 3, 4$ have zeroes at the points $(m+n\tau)\pi$, $(m+\frac{1}{2}+n\tau)\pi$, $(m+\frac{1}{2}+(n+\frac{1}{2})\tau)\pi$ and $(m+(n+\frac{1}{2})\tau)\pi$ respectively [12]. Taylor expansions of the theta – functions are the following [12]:

$$\theta_1(\pi z|\tau) = \pi z\theta_1'(0|\tau)\exp(-\sum_{j=1}^{\infty}\frac{1}{2j}\delta_{2j}z^{2j}),\ \theta_2(\pi z|\tau) = \theta_2(0|\tau)\exp(-\sum_{j=1}^{\infty}\frac{1}{2j}\alpha_{2j}z^{2j}),$$

$$\theta_3(\pi z|\tau) = \theta_3(0|\tau)\exp(-\sum_{j=1}^{\infty}\frac{1}{2j}\beta_{2j}z^{2j}),\ \text{and}\ \theta_4(\pi z|\tau) = \theta_4(0|\tau)\exp(-\sum_{j=1}^{\infty}\frac{1}{2j}\gamma_{2j}z^{2j})\ \text{(that is}$$

$$\theta_1(\pi z|\tau) = \pi z\theta_1'(0|\tau)(1 - \frac{\delta_2}{2}z^2 - (\frac{\delta_4}{4} - \frac{\delta_2^2}{8})z^4 + O(z^6)),\ \text{etc.), where:}$$

$$\delta_{2j}(\tau) = \sum_{\substack{n=-\infty\\|m|+|n|\neq 0}}^{\infty}\sum_{m=-\infty}^{\infty}\frac{1}{(m+n\tau)^{2j}} \qquad (13a),$$

$$\alpha_{2j}(\tau) = \sum_{n=-\infty}^{\infty} \sum_{m=-\infty}^{\infty} \frac{1}{(m+\frac{1}{2}+n\tau)^{2j}} \qquad (13b),$$

$$\beta_{2j}(\tau) = \sum_{n=-\infty}^{\infty} \sum_{m=-\infty}^{\infty} \frac{1}{(m+\frac{1}{2}+(n+\frac{1}{2})\tau)^{2j}} \qquad (13c)$$

$$\gamma_{2j}(\tau) = \sum_{n=-\infty}^{\infty} \sum_{m=-\infty}^{\infty} \frac{1}{(m+(n+\frac{1}{2})\tau)^{2j}} \qquad (13d).$$

The order of summation is important for these sums if $j=1$. (Of course, the Taylor expansions given above, albeit quite standard, are circular because they contain the values of $\theta_2(0,q)$, $\theta_3(0,q)$, $\theta_4(0,q)$ or $\theta_1'(0,q)$, which actually should be determined from the "real" Taylor expansions. Fortunately, this question is not relevant for the problems considered in the current Note).

From these formulae, we easily obtain the Taylor expansions of the logarithms of theta functions in some vicinity of zero, like

$$\ln[\theta_1(\pi z|\tau)/(\pi z \theta_1'(0|\tau))] = -\sum_{j=1}^{\infty} \frac{1}{2j} \delta_{2j} z^{2j} \qquad (14),$$

or $\ln[\theta_2(\pi z|\tau)/\theta_2(0|\tau)] = -\sum_{j=1}^{\infty} \frac{1}{2j} \alpha_{2j} z^{2j}$, etc. Note also useful for the future relation

$$\delta_2 = -\frac{\pi^2 \theta_1'''(0|\tau)}{3\theta_1'(0|\tau)} \qquad (15),$$

which immediately follows from Taylor series of $\theta_1(\pi z|\tau)$ given above.

Weierstrass sigma-function is defined as [12]

$$\sigma(z|\Lambda) = z \prod_{\omega \in \Lambda/\{0\}} (1-\frac{z}{\omega}) \exp(\frac{z}{\omega}+\frac{z^2}{2\omega^2}) \qquad (16),$$

where $\Lambda$ is a lattice formed by points $2m\omega_1 + 2n\omega_2$, with again $n, m \in Z$, and $\operatorname{Im}\frac{\omega_2}{\omega_1} > 0$. We have by definition $\zeta(z|\Lambda) = \frac{d}{dz}(\ln \sigma(z|\Lambda))$, $\wp(z|\Lambda) = -\frac{d\zeta(z|\Lambda)}{dz}$, $\wp_z(z|\Lambda) = \frac{d\wp(z|\Lambda)}{dz}$. Below we will use the notation $\sigma(z|\Lambda)$, $\zeta(z|\Lambda)$, etc. if an arbitrary lattice is using; if $\omega_1 = 1/2$, which is usually the case, we shall write $\sigma(z, \tau)$, $\zeta(z, \tau)$, etc.

The following Laurent expansion is well known (and evident from the above definition) [12]:

$$\zeta(z, \tau) = \frac{1}{z} - \sum_{k=2}^{\infty} \delta_{2k}(\tau) z^{2k-1} \qquad (17),$$

where $\delta_{2k}(\tau)$ are Eisenstein series defined by eq. (13a). To simplify the notation, below we may omit $\tau$ writing $\delta_{2k}(\tau)$ simply as $\delta_{2k}$. Thus $z\zeta(z, \tau) = 1 - \sum_{k=2}^{\infty} \delta_{2k} z^{2k}$ and $\ln(z\zeta(z, \tau)) = -\delta_4 z^4 - \delta_6 z^6 - (\delta_8 + \frac{1}{2}\delta_4^2) z^8 + O(z^{10})$, etc. Of course, we have

$$\wp(z, \tau) := -\frac{d}{dz}\zeta(z, \tau) = \frac{1}{z^2} + \sum_{k=2}^{\infty} (2k-1)\delta_{2k}(\tau) z^{2k-2} \qquad (18),$$

which is often written as

$$\wp(z, \tau) = \frac{1}{z^2} + z^2 \sum_{k=0}^{\infty} \frac{d_k}{k!} z^{2k} \qquad (18a),$$

where, evidently, $d_k = (2k+3)k!\delta_{2k+4}$.

Jacobi elliptical functions *sn(z, k), cn(z, k), dn(z, k)* and other similar, are scaled ratios of the appropriate theta-functions [12]. They will not be considered here.

*4.2 Equalities between different elliptic functions*

This is well known that all double-periodic meromorphic functions can be expressed via Weierstrass elliptical functions and their derivatives, and numerous concrete examples of the corresponding expressions/representations are commonly found in the literature [13 - 18]. By this reason, here we do not aim to consider many corresponding examples. Our purpose is to illustrate the principle and to concentrate on the cases which, by some reasons, seems are not widely indicated.

Our first example is the Weierstrass sigma function $\sigma(z|\Lambda)$. From the definition, it is clear that this is an entire function having simple zeroes at the lattice points. Thus the Weierstrass sigma function $\sigma(z, \tau)$ and the elliptic theta-function $\theta_1(\pi z|\tau)$ are both entire functions and have simple zeroes at $m+n\tau$. We know the Taylor expansion of the $\ln[\theta_1(\pi z|\tau)/(\pi z \theta_1'(0|\tau))]$, see eq. (14). Due to the Lemma 1, the Taylor expansion of the $\ln(\sigma(z, \tau)/z)$ must be exactly the same starting from the $O(z^4)$ term: $\ln(\sigma(z, \tau)/z) = a + bz^2 - \sum_{j=2}^{\infty} \frac{1}{2j} \delta_{2j} z^{2j}$. Directly from the definition (16) we see that $a=b=0$ hence

$$\sigma(z, \tau) = \exp\left(\frac{z^2}{2} \delta_2(\tau)\right) \frac{\theta_1(\pi z|\tau)}{\pi \theta_1'(0|\tau)} \qquad (19).$$

Scaling for an arbitrary $\sigma(z|\Lambda)$ immediately gives the "standard" formula [12]:

$$\sigma(z|\Lambda) = 2\omega_1 \exp\left(\frac{\eta_1 z^2}{2\omega_1}\right) \frac{\theta_1(\pi z/(2\omega_1), q)}{\pi \theta_1'(0, q)} \qquad (20),$$

where $\eta_1 = -\frac{\pi^2}{12\omega_1} \frac{\theta_1'''(0, q)}{\theta_1'(0, q)}$. Here and below, it might be useful to exploit [12]

$$\theta_1'(0, q) = 2q^{1/4} \prod_{j=1}^{\infty} (1-q^{2j})^3 \qquad (21).$$

**Remark 1.** *In light of the abovesaid, the Taylor expansion of the function $\frac{\sigma(z, \tau)}{z}$ is given by $\frac{\sigma(z, \tau)}{z} = \exp(-\frac{1}{4}\delta_4 z^4 - \frac{1}{6}\delta_6 z^6 - \frac{1}{8}\delta_8 z^8 + ...)$, sf. [12, 15]. To fully restore the known results, we need to add the formulae expressing $\delta_{2j}$ for j=4, 5, 6... via $\delta_4$ and $\delta_6$, which readily follow from e.g. the differential equation [12]*

$$\frac{d^2 \wp(z, \tau)}{dz^2} = 6\wp^2(z, \tau) - 30\delta_4 \qquad (22)$$

*by equating the Taylor series coefficients obtained by differentiating of eq. 18. Usually these formulae are written in the form including the coefficients $d_k$ of eq. (18a):*

$$\sum_{k=0}^{n} C_n^k d_k d_{n-k} = \frac{2n+9}{3n+6} d_{n+2} \qquad (23).$$

*(Differentiating eq. 18a, we see that the coefficient in front of the $z^{2n+4}$ term in the Taylor expansion of l.h.s of (22) is $\frac{(2n+5)(2n+6)d_{n+2}}{(n+2)!}$, while in the r.h.s it is*

$6 \sum_{k=0}^{n} \frac{d_k}{k!} \frac{d_{n-k}}{(n-k)!} + 12 \frac{d_{n+2}}{(n+2)!}$).

Our next example of the inverse application of the generalized Littlewood theorem to elliptic functions is the Weierstrass $\wp(z|\Lambda)$ functions. These functions have second order poles at the same points where simple zeroes of $\theta_1(\frac{\pi z}{2\omega_1}, q)$ function are located, which means that the function $\theta_1^{-2}(\frac{\pi z}{2\omega_1}, q)$ has the same poles as the Weierstrass $\wp(z|\Lambda)$ function. Here $\tau = \omega_2/\omega_1$ and $q = e^{i\pi\tau}$; again, we use this notation to be consistent with the data presented in [12].

But contrary to the entire theta function in the power of *-1*, the Weierstrass $\wp$-functions *have* zeroes. Situation is easiest with the functions $\wp(z|\Lambda) - e_i$, where $e_1 = \wp(\omega_1)$, $e_2 = \wp(\omega_2)$ and $e_3 = \wp(\omega_1 + \omega_2)$. "By construction", these functions have zeroes at points, coinciding with the "demi-lattice" points $m + \frac{1}{2} + n\tau$, $m + \frac{1}{2} + (n+\frac{1}{2})\tau$ or $m + (n+\frac{1}{2})\tau$ respectively, and it is clear that these are second order zeroes. This is also known that there are no other zeroes for functions $\wp(z|\Lambda) - e_i$ [12]. Thus, in particular, for the function $\wp(z|\Lambda) - e_1$, the positions of second order zeroes coincide with those of the entire function $\theta_2(\frac{\pi z}{2\omega_1}, q)$, and

Lemma 1 gives $\frac{d^2}{dz^2}\ln(\wp(z|\Lambda)-e_1) = \frac{d^2}{dz^2}\ln\frac{\theta_2^2(\pi z/(2\omega_1), q)}{\theta_1^2(\pi z/(2\omega_1), q)}$, so that

$\wp(z|\Lambda)-e_1 = C_2(\Lambda)\exp(C_1(\Lambda)z + C(\Lambda)z^2)\frac{\theta_2^2(\pi z/(2\omega_1), q)}{\theta_1^2(\pi z/(2\omega_1), q)}$. Coefficients $C(\Lambda)$ and $C_1(\Lambda)$ must be equal to zero for other coefficients are incompatible with the double periodicity of the both sides of the equality. (The ratio of the squares of theta-functions in the r.h.s. here is proportional to the square of Jacobi elliptic function $cs^2(\varsigma, k)$, see e.g. [12] for definitions). Analysis for the value of $z=0$ equating the coefficients in front of $1/z^2$ term in the Taylor expansions, shows $C_2 = \left(\frac{\pi \theta'^2_1(0,q)}{2\omega_1 \theta_2^2(0, q)}\right)^2$. Then we can use the known equality [12]:

$$\theta_1'(0, q) = \theta_2(0, q)\theta_3(0, q)\theta_4(0, q) \qquad (24)$$

to write the "standard" form:

$$\wp(z|\Lambda)-e_1 = \left(\frac{\pi\theta_3(0,q)\theta_4(0, q)}{2\omega_1}\right)^2 \frac{\theta_2^2(\pi z/(2\omega_1), q)}{\theta_1^2(\pi z/(2\omega_1), q)} \qquad (25).$$

Similar equalities can be in the same fashion established for other $\wp(z|\Lambda)-e_j$ functions, which we will not do here, see e.g. [12] for the list of formulae. We also will not exercise ourselves on the expression of the ratio of theta functions occurring in these formulae via appropriate Jacobi elliptic functions; see again [12]. Instead (see our motivation discussed above), we want to underline again that there is nothing special with the functions $\wp(z|\Lambda)-e_j$. The following Theorem, which seems is almost never appears in the monographs devoted to the elliptic functions, holds true.

**Theorem 2.** *Let* a *be an arbitrary complex number not equal to* $e_1, e_2, e_3$ *defined above, and numbers* $\alpha_{1,2}$ *be such that the equality* $\wp(\alpha_j | \Lambda) = a$, *j=1,2, holds true. Let also the difference* $\alpha_1 - \alpha_2 \neq 2n\omega_1 + 2m\omega_2$. *Then*

$$\wp(z|\Lambda) - a = C \frac{\theta_1(\pi z/(2\omega_1) - \pi\alpha_1/(2\omega_1))\theta_1(\pi z/(2\omega_1) - \pi\alpha_2/(2\omega_1)),\ q)}{\theta_1^2(\pi z/(2\omega_1),\ q)} \quad (26),$$

*where* $\tau = \omega_2/\omega_1$, $q = e^{i\pi\tau}$,

$$C = \frac{\pi^2}{4\omega_1^2} \frac{\theta_1'^2(0,\ q)}{\theta_1(\pi\alpha_1/(2\omega_1),\ q)\theta_1(\pi\alpha_2/(2\omega_1),\ q)} \quad (27).$$

**Proof.** First, we note that due to the argument principle, double-periodic function $\wp(z|\Lambda) - a$, having in each fundamental parallelogram one pole of the second order, has in this parallelogram exactly two zeroes taking into account their multiplicity. (Exactly for $a=e_j$, $j=1, 2, 3$, excluded by the Theorem conditions, this function has one double zero inside the fundamental parallelogram). Thus for any $a$, we can find the numbers $\alpha_{1,2}$ requested by the Theorem. Clearly, the functions $\wp(z|\Lambda) - a$ and $\theta_1^{-2}(\pi z/(2\omega_1),\ q)$ have the same poles of the second order at the lattice points. Function $\wp(z|\Lambda) - a$ additionally has simple zeroes at the points $\alpha_{1,2} + 2n\omega_1 + 2m\omega_2$, and the product $\theta_1(\pi z/(2\omega_1) - \pi\alpha_1/(2\omega_1),\ q) \times \theta_1(\pi z/(2\omega_1) - \pi\alpha_2/(2\omega_1),\ q)$ is an entire function also having simple zeroes exactly at these points. The functions involved are both double-periodic hence the Lemma 1 is applicable already for *n*=2. The double integration of the obtained equality of the second derivatives of the corresponding logarithms gives

$$\wp(z|\Lambda) - a = C(\omega_1,\ q) \frac{\theta_1(\pi z/(2\omega_1) - \pi\alpha_1/(2\omega_1),\ q)\theta_1(\pi z/(2\omega_1) - \pi\alpha_2/(2\omega_1)),\ q)}{\theta_1^2(\pi z/(2\omega_1),\ q)}.$$ (Again,

the possible factor $\exp(Bz) \equiv 1$ due to the double-periodicity of the functions involved). The coefficient *C* is established by analyzing the case when *z* tends to zero by equating the terms proportional to $1/z^2$ in Laurent expansions. We have
$$C = \frac{\pi^2}{4\omega_1^2} \frac{\theta_1'^2(0,\, q)}{\theta_1(\pi\alpha_1/(2\omega_1),\, q)\theta_1(\pi\alpha_2/(2\omega_1),\, q)},$$
and the theorem is proven. □

**Corollary 1.** *Applying (26), (27) for a=0, we have the following presentation of the Werierstrass $\wp(z)$ function via $\theta_1$-elliptical functions:*

$$\wp(z\,|\,\Lambda) = \frac{\pi^2}{4\omega_1^2} \frac{\theta_1'^2(0,\, q)}{\theta_1(\pi\alpha_1/(2\omega_1),\, q)\theta_1(\pi\alpha_2/(2\omega_1),\, q)} \times$$

$$\frac{\theta_1(\pi z/(2\omega_1) - \pi\alpha_1/(2\omega_1),\, q)\theta_1(\pi z/(2\omega_1) - \pi\alpha_2/(2\omega_1)),\, q)}{\theta_1^2(\pi z/(2\omega_1),\, q)} \quad (28),$$

*where $\alpha_{1,2}$ are some solutions of the equation $\wp(z\,|\,\Lambda) = 0$ such that $\alpha_1 - \alpha_2 \neq 2n\omega_1 + 2m\omega_2$. In particular, we can take these solutions as discussed (for $\omega_1 = 1/2$) in Refs. [19, 20] to be equal $\pm\alpha_1$; here $\alpha_1$ belongs to the first fundamental parallelogram . Then*

$$\wp(z\,|\,\Lambda) = -\pi^2 \frac{\theta_1'^2(0,\, q)}{\theta_1^2(\pi\alpha_1,\, q)} \times \frac{\theta_1(\pi z - \pi\alpha_1,\, q)\theta_1(\pi z + \pi\alpha_1,\, q)}{\theta_1^2(\pi z,\, q)}.$$

**Remark 2**. *It is instructive to see what happens in the limit $\tau \to i\infty$ when $\wp(z,\, \tau)$ tends to $\wp(z) = \frac{\pi^2}{\sin^2 \pi z} - \frac{\pi^2}{3}$. The smallest in the module solutions of $\wp(z)=0$ are then given by the complex numbers $\alpha_{1,2} = \pm(\frac{1}{2} + \frac{i}{2\pi}\ln(5 + 2\sqrt{6}))$ - the values still not easy to work with. Thus for even greater simplicity, let us consider more transparent case e.g. of $\wp(z) - \pi^2 = \frac{\pi^2}{\sin^2 \pi z} - \frac{4\pi^2}{3}$. We have:*

$$\wp(z) - \pi^2 = \frac{\pi^2}{\sin^2 \pi z} - \frac{4\pi^2}{3} = -\frac{4\pi^2}{3\sin^2 \pi z}\sin(\pi z - \frac{\pi}{3})\sin(\pi z + \frac{\pi}{3}),$$ *and this gives the aforementioned presentation for the case.*

The derivative of the Weierstrass $\wp$-function $\wp_z(z, \tau)$ has third order poles at $k + j\tau$ and zeroes at $k + 1/2 + (j + 1/2)\tau$ [12]. Thus the lattice, formed by the irregular points of the logarithm (i.e. poles *and* zeroes) of the $\wp_z(z, \tau)$ coincides with the zeroes of the function $\sigma(2z, \tau)$. We can convert the third order poles of the $\wp_z(z, \tau)$ function into simple zeroes by multiplying it by $\sigma^4(z, \tau)$ function: thus both functions $\sigma(2z, \tau)$ and $\wp_z(z, \tau)\sigma^4(z, \tau)$ are entire and have the same simple zeroes. Using our method, we than can easily establish the well-known equality

$$\sigma(2z, \tau) = -\wp_z(z, \tau)\sigma^4(z, \tau) \qquad (29).$$

Of course, there are also infinitely many other possibilities of such type. For instance, we see that the entire function $\sigma(z + \frac{1}{2}, \tau)$ has simple zeroes at $m + \frac{1}{2} + n\tau$ points, thus at the same points as the entire function $\theta_2(\pi z | \tau)$, so that $\sigma(z + \frac{1}{2}, \tau) = C\exp(C_1 z + C_2 z^2)\theta_2(\pi z | \tau)$. Omitting clear from the abovesaid details, we arrive to

$$\sigma(z + \frac{1}{2}, \tau) = \frac{\sigma(1/2, \tau)}{\theta_2(0 | \tau)} \exp(\zeta(1/2, \tau)z + \frac{1}{2}\delta_2(\tau)z^2)\theta_2(\pi z | \tau) \qquad (30).$$

Similar equalities hold for the functions $\sigma(z + \frac{\tau}{2}, \tau)$, $\sigma(z + \frac{1}{2} + \frac{\tau}{2}, \tau)$, see e.g. [12].

Next example seems is apparently also not quite common. We see that $\sigma(2z, \tau)$ and $\theta_1(\pi z, \tau)\theta_2(\pi z, \tau)\theta_3(\pi z, \tau)\theta_4(\pi z, \tau)$ are both entire and have the same simple zeroes. Thus we have

$$8\frac{d^3}{dz^3}\ln(\sigma(2z, \tau)) = \pi^3 \frac{d^3}{dz^3}\ln[\theta_1(\pi z, \tau)\theta_2(\pi z, \tau)\theta_3(\pi z, \tau)\theta_4(\pi z, \tau)] \text{ and}$$

$4\dfrac{d^2}{dz^2}\ln(\sigma(2z,\ \tau))=\pi^2\dfrac{d^2}{dz^2}\ln[\theta_1(\pi z,\ \tau)\theta_2(\pi z,\ \tau)\theta_3(\pi z,\ \tau)\theta_4(\pi z,\ \tau)]+C_1$. Equating when $z$ tends to zero, we get $0=C_1-\delta_2-\alpha_2-\beta_2-\gamma_2$, see eqs. (13) for definitions. In Ref. [9], we showed that $3\delta_2=\alpha_2+\beta_2+\gamma_2$, thus $C_1=4\delta_2$. For the first derivatives

$2\dfrac{d}{dz}\ln(\sigma(2z,\ \tau))=\pi\dfrac{d}{dz}\ln[\theta_1(\pi z,\ \tau)\theta_2(\pi z,\ \tau)\theta_3(\pi z,\ \tau)\theta_4(\pi z,\ \tau)]+C_1 z+C_2$, so that trivially $C_2=0$. Thus $\sigma(2z,\ \tau)=C\exp(2\delta_2 z^2)\theta_1(\pi z,\ \tau)\theta_2(\pi z,\ \tau)\theta_3(\pi z,\ \tau)\theta_4(\pi z,\ \tau)$. Equating functions at zero, we have $2=C\pi\theta_1'(0|\tau)\theta_2(0|\tau)\theta_3(0|\tau)\theta_4(0|\tau)$ and with eq. (25), $C=2/(\pi[\theta_2(0)\theta_3(0)\theta_4(0)]^2)$, hence

$$\sigma(2z,\ \tau)=\dfrac{2}{\pi\theta_2^2(0)\theta_3^2(0)\theta_4^2(0)}\exp(2\delta_2 z^2)\theta_1(\pi z,\ \tau)\theta_2(\pi z,\ \tau)\theta_3(\pi z,\ \tau)\theta_4(\pi z,\ \tau) \qquad (31).$$

As usual, $\delta_2$ can be expressed via the derivatives ratio using eq. (15). Now note that by eq. (19) $\sigma(2z,\ \tau)=\exp\!\left(2z^2\delta_2(\tau)\right)\dfrac{\theta_1(2\pi z|\tau)}{\pi\theta_1'(0|\tau)}$, so that we have the following expression

$$\theta_1(2z|\tau)=\dfrac{2}{\theta_2(0|\tau)\theta_3(0|\tau)\theta_4(0|\tau)}\theta_1(z|\tau)\theta_2(z|\tau)\theta_3(z|\tau)\theta_4(z|\tau) \qquad (32)$$

- which, of course, can be easily established by our approach also in the "direct way". See the discussion of the transformation rules for the elliptical functions in the next sub-section.

Analogously, and this is our final example in this sub-section, we see that $\wp_z(z,\ \tau)$ and $\theta_1^{-3}(\pi z|\tau)\theta_2(\pi z|\tau)\theta_3(\pi z|\tau)\theta_4(\pi z|\tau)$ have the same poles and zeroes. Omitting standard easy details, we get

$\wp_z(z,\ \tau)=C\theta_1^{-3}(\pi z|\tau)\theta_2(\pi z|\tau)\theta_3(\pi z|\tau)\theta_4(\pi z|\tau)$, and equating at zero, we have

$-2=C\pi^{-3}[\theta_1'(0|\tau)]^{-3}\theta_2(0|\tau)\theta_3(0|\tau)\theta_4(0|\tau)=C\pi^{-3}[\theta_2(0|\tau)\theta_3(0|\tau)\theta_4(0|\tau)]^{-2}$.

Hence $C=-2\pi^3[\theta_2(0|\tau)\theta_3(0|\tau)\theta_4(0|\tau)]^2$ and

$$\wp_z(z) = -\frac{2\pi^3}{[\theta_2(0)\theta_3(0)\theta_4(0)]^2}\theta_1^{-3}(\pi z|\tau)\theta_2(\pi z, \tau)\theta_3(\pi z, \tau)\theta_4(\pi z, \tau) \quad (33).$$

Here in the r.h.s. the expressions of $\frac{\theta_2(\pi z|\tau)}{\theta_1(\pi z|\tau)}$, $\frac{\theta_3(\pi z|\tau)}{\theta_1(\pi z|\tau)}$ and $\frac{\theta_4(\pi z|\tau)}{\theta_1(\pi z|\tau)}$ via Jacobi elliptical functions $cs(\varsigma, k)$, $ds(\varsigma, k)$, $ns(\varsigma, k)$ respectively, see [12] for details, can be used.

### 4.3. n-tuple products

Let us consider entire theta function $\theta_1(nz|n\tau)$ with the integer $n=2, 3, 4...$ It has zeroes at the points $nz = k\pi + jn\pi\tau$, i.e. $z = \frac{k}{n}\pi + j\pi\tau$ (as usual, $k, j \in Z$), and thus evidently has the same zeroes as the product $\prod_{k=0}^{n-1}\theta_1(z + \frac{k\pi}{n}|\tau)$. This, together with the asymptotic of the functions involved for large values of $|z|$, immediately enable to write, by applying Lemma 1:

$$\theta_1(nz|n\tau) = C(n, \tau)\exp(C_1(n,\tau)z + C_2(n,\tau)z^2)\prod_{k=0}^{n-1}\theta_1(z + \frac{k\pi}{n}|\tau).$$ $C_1(n, \tau)$ and $C_2(n, \tau)$ are clearly zero because other values are incompatible with the $z$-periodicity of the functions involved. (Purely imaginary $C(n, \tau) = ik$ with integer $k$ might be compatible (for we do not have *double* periodicity now), but such value cannot be obtained by comparison of the real Taylor expansions at $z = 0$). The constant $C(n, \tau)$ should be determined from this equality applied at any particular value of $z$. Using for $z$ tending to zero $\theta_1(nz|n\tau) = nz\theta'_1(0|n\tau) + O(z^2)$, we have $nz\theta'_1(0|n\tau) = C(n,\tau)\theta'_1(0|\tau)z\prod_{k=1}^{n-1}\theta_1(\frac{k\pi}{n}|\tau)$ hence

$$\theta_1(nz|n\tau) = \frac{n\theta'_1(0|n\tau)}{\theta'_1(0|\tau)\prod_{k=1}^{n-1}\theta_1(\frac{k\pi}{n}|\tau)}\prod_{k=0}^{n-1}\theta_1(z + \frac{k\pi}{n}|\tau) \quad (34).$$

In particular, $\theta_1(2z|2\tau) = \dfrac{2\theta_1'(0|2\tau)\theta_1(z|\tau)\theta_1(z+\frac{\pi}{2}|\tau)}{\theta_1'(0|\tau)\theta_1(\frac{\pi}{2}|\tau)} = \dfrac{2\theta_1'(0|2\tau)\theta_1(z|\tau)\theta_2(z|\tau)}{\theta_1'(0|\tau)\theta_2(0|\tau)}$ and

$$\theta_1(4z|4\tau) = \dfrac{4\theta_1'(0|4\tau)\theta_1(z|\tau)\theta_1(z+\frac{\pi}{4}|\tau)\theta_1(z+\frac{\pi}{2}|\tau)\theta_1(z+\frac{3\pi}{4}|\tau)}{\theta_1'(0|\tau)\theta_1(\frac{\pi}{4}|\tau)\theta_1(\frac{\pi}{2}|\tau)\theta_1(\frac{3\pi}{4}|\tau)}$$

$$= \dfrac{4\theta_1'(0|4\tau)\theta_1(z|\tau)\theta_1(z+\frac{\pi}{4}|\tau)\theta_2(z|\tau)\theta_1(z+\frac{3\pi}{4}|\tau)}{\theta_1'(0|\tau)\theta_1(\frac{\pi}{4}|\tau)\theta_2(0|\tau)\theta_1(\frac{3\pi}{4}|\tau)}.$$

**Remark 3.** *By unclear reasons this relation is almost never presented in such a (general) form and, if presented for some particular case, is given in somewhat "artificial" form - like, for example, the following Landen transformations [12]:*

$\theta_1(4z|4\tau) = \dfrac{\theta_1(z|\tau)\theta_1(\frac{\pi}{4}-z|\tau)\theta_1(\frac{\pi}{4}+z|\tau)\theta_2(z|\tau)}{\theta_3(0|\tau)\theta_4(0|\tau)\theta_3(\frac{\pi}{4}|\tau)}$, *or* $\theta_1(2z|2\tau) = \dfrac{\theta_1(z|\tau)\theta_2(z|\tau)}{\theta_4(0|2\tau)}$. *(This is understandable, because Landen transformations historically came from the manipulations with the elliptical integrals achieved mainly with variable changes [10, 12 - 17], but still). Of course, the equivalence of the transformations can be shown using formulae connecting the values of the elliptical functions and their derivatives at certain values, see [12 -17], especially [16], but this is actually not easy. Alternatively, the comparison of these formulae can be seen as a* proof *of certain relations like, say,* $\dfrac{2\theta_1'(0|2\tau)}{\theta_1'(0|\tau)\theta_1(\frac{\pi}{2}|\tau)} = \dfrac{1}{\theta_4(0|2\tau)}$.

*Further, let us note the following.* $\lim_{\tau\to i\infty} \theta_1(z|\tau) = 0$ *for any finite z. But from the definition it is easy to see that* $\lim_{\tau\to i\infty} e^{-i\pi\tau/4}\theta_1(z|\tau) = 2\sin(z)$. *Applying this (properly scaled) limit to the both sides of (34), we get eq. (8).*

There is also quite different possibility. We can consider not $\theta_1(nz|n\tau)$ but $\theta_1(nz|\tau)$ instead and obtain the following Theorem.

**Theorem 3.** *For odd 2l+1=3, 5, 7..., we have:*

$$\theta_1((2l+1)z|\tau) = C\prod_{k=-l}^{l}\prod_{j=-l}^{l}\theta_1(z+\frac{k\pi}{2l+1}+\frac{j\pi\tau}{2l+1}|\tau) \qquad (35)$$

with

$$C^{-1} = \frac{1}{2l+1}\prod_{k=-l}^{l}\prod_{\substack{j=-l,\\|j|+|l|\neq 0}}^{l}\theta_1(\frac{k\pi}{2l+1}+\frac{j\pi\tau}{2l+1}|\tau) \qquad (36).$$

*For even 2l=2, 4, 6..., we have*

$$\theta_1(2lz,\ \tau) = \tilde{C}\prod_{m=-(l-1)}^{l}\theta_4(z+\frac{m\pi}{2l})\times\prod_{k=-(l-1)}^{l}\prod_{j=-(l-1)}^{l-1}\theta_1(z+\frac{k\pi}{2l}+\frac{j\pi\tau}{2l}|\tau) \qquad (37)$$

with

$$\tilde{C}^{-1} = \frac{1}{2l}\prod_{m=-(l-1)}^{l}\theta_4(\frac{m\pi}{2l})\times\prod_{k=-(l-1)}^{l}\prod_{\substack{j=-(l-1),\\|j|+|k|\neq 0}}^{l-1}\theta_1(\frac{k\pi}{2l}+\frac{j\pi\tau}{2l}|\tau) \qquad (38).$$

**Proof.** This is easy to see that the functions $\theta_1(nz|\tau)$ and $\prod_{k=0}^{n-1}\prod_{j=0}^{n-1}\theta_1(z+\frac{k\pi}{n}+\frac{j\pi\tau}{n}|\tau)$ are both entire and have the same zeroes, but to find an exponential factor we need to "symmetrize" this product.

Let first $n$ be odd, $n=2l+1$. Than $\theta_1((2l+1)z|\tau)$ has the same zeroes and (no) poles as the "symmetric" product $\prod_{k=-l}^{l}\prod_{j=-l}^{l}\theta_1(z+\frac{k\pi}{n}+\frac{j\pi\tau}{n}|\tau)$. Thus by the Lemma 1 we get

$$\theta_1((2l+1)z|\tau) = C\exp(C_1z+C_2z^2)\prod_{k=-l}^{l}\prod_{j=-l}^{l}\theta_1(z+\frac{k\pi}{n}+\frac{j\pi\tau}{n}|\tau) \qquad (39).$$

As usual, the coefficient $C_2$ should be equal to zero due to the $z$-periodicity. The coefficient $C_1$ is equal to zero due to the even nature of the $\theta_1(nz|\tau)/z$ and $\frac{1}{z}\prod_{k=-l}^{l}\prod_{j=-l}^{l}\theta_1(z+\frac{k\pi}{n}+\frac{j\pi\tau}{n}|\tau)$ functions. This is seen as follows: $\theta_1(nz|\tau)/z$ is even, and all factors in the $\prod_{k=-l}^{l}\prod_{j=-l}^{l}\theta_1(z+\frac{k\pi}{n}+\frac{j\pi\tau}{n}|\tau)$ come in pairs with $\pm k$, $\pm j$ if $k\neq 0$, $j\neq 0$. Namely, we have

$$\theta_1(z+\frac{k\pi}{n}+\frac{j\pi}{n}|\tau)\theta_1(z-\frac{k\pi}{n}-\frac{j\pi}{n}) = -\theta_1(z+\frac{k\pi}{n}+\frac{j\pi}{n}|\tau)\theta_1(-z+\frac{k\pi}{n}+\frac{j\pi}{n})$$

so that this product is an even function. If $k$ or $j$ are equal to zero, this is also valid. If both $k=j=0$ simultaneously, the factor $\theta_1(z|\tau)/z$ is even. Thus, comparing the $O(z)$ terms of the Taylor development of the logarithms of both sides of (39), we get $C_1=0$. Coefficient $C$ is obtained equating the values at zero:

$$C^{-1} = \frac{1}{2l+1}\prod_{k=-l}^{l}\prod_{\substack{j=-l,\\|j|+|l|\neq 0}}^{l}\theta_1(\frac{k\pi}{n}+\frac{j\pi\tau}{n}|\tau).$$

What does happen if $n=2l$ is even? Then in passing from $\prod_{k=0}^{2l-1}\prod_{j=0}^{2l-1}\theta_1(z+\frac{k\pi}{n}+\frac{j\pi\tau}{n}|\tau)$ to the symmetrized form, we obtain still not the fully symmetric product $\prod_{k=-(l-1)}^{l}\prod_{j=-(l-1)}^{l}\theta_1(z+\frac{k\pi}{n}+\frac{j\pi\tau}{n}|\tau)$. The situation with $\frac{1}{z}\prod_{k=-(l-1)}^{l-1}\prod_{j=-(l-1)}^{l-1}\theta_1(z+\frac{k\pi}{n}+\frac{j\pi\tau}{n}|\tau)$ is the same as before: this is an even function, but now we get also the "unpaired" factors at $k=l$ and $j=l$. These are the functions $\theta_1(z+\frac{\pi}{2}+\frac{\pi\tau j}{n}|\tau)$ and $\theta_1(z+\frac{\pi k}{n}+\frac{\pi\tau}{2}|\tau)$ respectively. They are equal respectively to $\theta_2(z+\frac{\pi\tau j}{n}|\tau)$ and $Az^{-iz}\theta_4(z+\frac{\pi k}{n}|\tau)$, where the constant $A$ easily follows from [12]

$$\theta_1(z+\frac{\pi\tau}{2}|\tau) = -ie^{-iz}e^{-i\pi\tau/4}\theta_4(z|\tau) \tag{40},$$

but for our current purposes important is only that it does not depend on $z$. The functions $\theta_2(z+\frac{\pi\tau j}{n}|\tau)$ with $1\le j\le l-1$, when paired as

$\theta_2(z+\frac{\pi\tau j}{n}|\tau)\theta_2(z-\frac{\pi\tau j}{n})=\theta_2(z+\frac{\pi\tau j}{n}|\tau)\theta_2(-z+\frac{\pi\tau j}{n})$, again are even functions. The function $\theta_2(z|\tau)$, corresponding to $j=0$, is also even.

Thus among the functions corresponding to $k=l$, only $\theta_1(z+\frac{\pi}{2}+\frac{\pi\tau}{2}|\tau)$ rests unpaired. It is equal to $Be^{-iz}\theta_3(z|\tau)$ [12], where $B$ does not depend on $z$; the function $\theta_3(z|\tau)$ is even.

In the second group of functions (corresponding to $j=l$), the functions $Az^{-iz}\theta_4(z+\frac{\pi k}{n}|\tau)$ for $k\ne 0$ are also paired:

$A_1e^{-iz}\theta_4(z+\frac{\pi k}{n}|\tau)A_2e^{-iz}\theta_4(z-\frac{\pi k}{2})=A_1A_2e^{-2iz}\theta_4(z+\frac{\pi k}{n}|\tau)\theta_4(-z+\frac{\pi k}{2})$. For $k=0$, we have $Ae^{-iz}\theta_4(z|\tau)$ with an even function $\theta_4(z|\tau)$; the case $k=j=l$ has been already considered just before.

Thus we see that the product of all these unpaired factors, related with the theta functions at $k=l$ and $j=l$, is $e^{-2liz}\varphi(z)$, where $\varphi(z)$ is an even function. Thus comparing the $O(z)$ terms of Taylor expansions of the logarithms of the both sides of the equation $\theta_1(2lz|\tau)=C\exp(C_1z)\prod_{k=-(l-1)}^{l}\prod_{j=-(l-1)}^{l}\theta_1(z+\frac{k\pi}{n}+\frac{j\pi\tau}{n}|\tau)$, we obtain $C_1=2il$.

We might finish the consideration at this point, but it seems preferable to give the purely real form using (40) again:

$$\theta_1(2lz, \tau) = \tilde{C} \prod_{m=-(l-1)}^{l} \theta_4(z + \frac{m\pi}{2l}) \times \prod_{k=-(l-1)}^{l} \prod_{j=-(l-1)}^{l-1} \theta_1(z + \frac{k\pi}{2l} + \frac{j\pi\tau}{2l} | \tau).$$ (Note, that for $m=l$ in the first product, $\theta_4(z + \frac{\pi}{2}) = \theta_3(z)$ [12], we used this earlier). Coefficient $C$ is obtained equating the values of both sides of the above equation at zero:

$$\tilde{C}^{-1} = \frac{1}{2l} \prod_{m=-(l-1)}^{l} \theta_4(\frac{m\pi}{2l}) \times \prod_{k=-(l-1)}^{l} \prod_{\substack{j=-(l-1),\\|j|+|k|\neq 0}}^{l-1} \theta_1(\frac{k\pi}{2l} + \frac{j\pi\tau}{2l} | \tau). \quad \square$$

**Remark 4.** 1. *The appearance of the factor $e^{2ilz}$ during the derivation of the n-tuple relation for even $n=2l$ might look unexpected. However, we get it already for $n=2$. Directly from the Lemma 1, we easily prove* $\theta_1(2z | \tau) = \frac{2\theta_1(z | \tau)\theta_2(z | \tau)\theta_3(z | \tau)\theta_4(z | \tau)}{\theta_2(0 | \tau)\theta_3(0 | \tau)\theta_4(0 | \tau)}.$
*From our Theorem 3 we have*

$$\theta_1(2z | \tau) = Ce^{2iz}\theta_1(z | \tau)\theta_1(z + \frac{\pi}{2} | \tau)\theta_1(z + \frac{\pi\tau}{2} | \tau)\theta_1(z + \frac{\pi}{2} + \frac{\pi\tau}{2} | \tau)$$ - *and this is consistent, because* $\theta_1(z + \frac{\pi}{2} | \tau) = \theta_2(z)$, $\theta_1(z + \frac{\pi\tau}{2} | \tau) = Ae^{-iz}\theta_4(z | \tau)$, $\theta_1(z + \frac{\pi}{2} + \frac{\pi\tau}{2} | \tau) = Be^{-iz}\theta_3(z | \tau).$

*2. Along the same lines, relations similar to (35) – (38) can be established for the Weierstrass sigma functions $\sigma(nz, n\tau)$, $\sigma(nz, \tau)$. These same relations also can be obtained from the above relations for the first theta-function, supplemented with eq. (19) or (20). The present author, however, was able to find only the relation analogous to (35) for the Weierstrass $\sigma(nz, \tau)$ function, formulae 23.10.13 - 23.10.16 in [12], (with a misprint in 23.10.14: instead of*
$A_n = n \prod_{j=0}^{n-1} \prod_{\substack{l=0\\|j|+|l|\neq 0}}^{n-1} \sigma^{-1}(\frac{2j\omega_1}{n} + \frac{2l\omega_3}{n})$ *it is written* $A_n = n \prod_{j=0}^{n-1} \prod_{\substack{l=0\\l\neq j}}^{n-1} \sigma^{-1}(\frac{2j\omega_1}{n} + \frac{2l\omega_3}{n})$), *but not for the theta-functions.*

To finish this Section, let us briefly present an easy case of the $n$-tuple relations for the $\wp_z$ function. We have the following simple Theorem.

**Theorem 4.** *For n=2, 3, 4... we have*

$$\wp_z(nz, n\tau) = \frac{1}{n^3 \prod_{k=1}^{n-1} \wp_z(\frac{k\pi}{n}, \tau)} \prod_{k=0}^{n-1} \wp_z(z + \frac{k}{n}, \tau) \tag{41}$$

*and*

$$\wp_z(nz \mid \tau) = C \prod_{k=0}^{n-1} \prod_{j=0}^{n-1} \wp_z(z + \frac{k}{n} + \frac{j\tau}{n}, \tau) \tag{42}$$

*with*

$$C^{-1} = n^3 \prod_{k=0}^{n-1} \prod_{\substack{j=0, \\ |j|+|l| \neq 0}}^{n-1} \wp_z(\frac{k}{n} + \frac{j\tau}{n}, \tau) \tag{43}.$$

**Proof.** Above, we have already noticed that the function $\wp_z(z, \tau)$ has third order poles at $z = k + j\tau$ and simple zeroes at $z = k + 1/2 + (j + 1/2)\tau$ [12]. Such a regular lattice of peculiarities guaranties that the products in both sides of above equations have the same poles and zeroes. Application of Lemma 1 immediately gives

$$\wp_z(nz, n\tau) = C \exp(C_1 z) \prod_{k=0}^{n-1} \wp_z(z + \frac{k}{n}, \tau) \text{ and}$$

$\wp_z(nz \mid \tau) = \tilde{C} \exp(\tilde{C}_1 z) \prod_{k=0}^{n} \prod_{j=0}^{n-1} \wp_z(z + \frac{k}{n} + \frac{j\tau}{n}, \tau)$. (As during the proof of the Theorem 2, the function $\wp_z(z, \tau)$ is double-periodic hence already the integrals $\int_{C_i} \frac{1}{z^2} \ln(\wp_z(z, \tau)) dz$ tend to zero in the limit of infinitely large contours; starting with them, we do not obtain the terms $\exp(C_2 z^2)$ at all). The same double periodicity of the involved functions ensures that $C_1$, $\tilde{C}_1 = 0$. Equating the values of the both sides of these equations at z=0 (remind that $\wp_z(z, \tau) = \frac{d}{dz}\wp(z, \tau) = -\frac{2}{z^3} + O(z)$), we get (41) - (43). □

*4.4. Fundamental modular transformations and Jacobi's triple product*

Let us briefly discuss the fundamental modular transformations. Analogously to what has been done above, we can easily establish $\theta_1(z|\tau+1) = C_1\theta_1(z|\tau)$: the possible exponential factor $\exp(Az^2 + Bz) \equiv 1$. (Any $A \neq 0$ is incompatible with the z-periodicity, while $B=0$ simply due to the evenness of the functions $\ln(\theta_1(z|\tau)/z)$).

Similarly, $\theta_1(z|-\frac{1}{\tau}) = C_2 \exp(\frac{i\tau z^2}{\pi})\theta_1(\tau z|\tau)$. Here in the possible factor $\exp(Az^2 + Bz)$, $B=0$ for the same reason as for the previous relation, while to establish the coefficient $A$ we used the Taylor expansion (13a). But this calculation requires certain caution! Clearly,

$$\delta_2(-1/\tau) = \sum_{n=-\infty}^{\infty}\sum_{\substack{m=-\infty \\ |m|+|n|\neq 0}}^{\infty} \frac{1}{(m-n/\tau)^2} = \tau^2 \sum_{n=-\infty}^{\infty}\sum_{\substack{m=-\infty \\ |m|+|n|\neq 0}}^{\infty} \frac{1}{(m\tau+n)^2} \qquad (44),$$

where $\delta_2(\tau)$ is defined by (13a), so that one might decide that $\frac{d^2}{dz^2}\ln(\theta_1(z|-\frac{1}{\tau})/z)|_{z=0} = \frac{d^2}{dz^2}\ln(\theta_1(\tau z|\tau)/z)|_{z=0}$ and $A=0$ which is clearly impossible: $\theta_1(z|-\frac{1}{\tau})$ and $\theta_1(\tau z|\tau)$ have different incompatible periods. However, the *order* of the summation in the r.h.s of eq. (39) is different from that of eq. (13a), so we need to remind the Eisenstein relation $\sum_{m=-\infty}^{\infty}\sum_{\substack{n=-\infty \\ |m|+|n|\neq 0}}^{\infty} \frac{1}{(m+n\tau)^2} = \delta_2(\tau) - \frac{2\pi i}{\tau}$ (see the discussion and further references in [9]), and in such a way we obtain $A = \frac{i\tau}{\pi}$.

Of course, $C_1 = \frac{\theta_1'(0|1+\tau)}{\theta_1'(0|\tau)}$ and $C_2 = \frac{\theta_1'(0|-1/\tau)}{\tau\theta_1'(0|\tau)}$. To move further, the use of the original definitions eqs. (11) is necessary to get [12]:

$$\theta_1(z|1+\tau) = i^{1/2}\theta_1(z|\tau) \qquad (45),$$

$$\theta_1(z|-\frac{1}{\tau}) = -i(-i\tau)^{1/2} \exp(\frac{i\tau z^2}{\pi})\theta_1(\tau z|\tau) \qquad (46).$$

Finally, let us look at Jacobi triple product. It is not too difficult to note that a whole function, having zeroes at $m+j\tau$, that is coinciding with those of the function $\theta_1(\pi z|\tau)$, can be written as $\psi_1(z,\tau) = \prod_{j=0}^{\infty}(1-e^{2\pi i((j+1)\tau-z)})(1-e^{2\pi i(j\tau+z)})$, and similar expressions can be easily obtained for other theta-functions. (This might look not too familiar, because more standard terminology/notation in the field is the use of $x = e^{2\pi i z}$, $q = e^{2\pi i\tau}$ (or $x = e^{\pi i z}$ and the nome $q = e^{\pi i\tau}$) and $(x;q) = \prod_{j=0}^{\infty}(1-xq^j)$ - then we see that $\psi_1(z,\tau) \equiv (x;q)(q/x;q)$, etc.). Application of Lemma 1 immediately enables to write (in somewhat weird notation):

$\theta_1(\pi z|\tau) = C(\tau)\exp(C_1(\tau)z)(x;q)(q/x;q)$, where possible quadratic $\exp(C_2(\tau)z^2)$ term is again excluded by periodicity. Exact value of $C_1$ is calculated comparing the values of the derivatives of the logarithms at $z=0$. Trivially, $\frac{d}{dz}(\ln\frac{\theta_1(\pi z|\tau)}{z})|_{z=0} = 0$, while for the r.h.s we have $(x;q)(q/x;q)/z = \prod_{j=0}^{\infty}(1-xq^j)(1-q^{j+1}/x)/z = \frac{1-x}{z}\prod_{j=1}^{\infty}(1-xq^j)(1-q^j/x)$,

and further $\lim_{x\to 1}\frac{d}{dx}\ln(\frac{1-x}{z}) = \lim_{x\to 1}\frac{d}{dx}\ln\frac{1-x}{2\pi i\ln x} = \lim_{\delta\to 0}\frac{d}{d\delta}\ln\frac{-\delta}{\ln(1+\delta)} = \frac{1}{2}$, together

with $\lim_{x\to 1}(\frac{d}{dx}\ln(\prod_{j=1}^{\infty}(1-xq^j)(1-q^j/x))) = \lim_{x\to 1}(\sum_{j=1}^{\infty}(-\frac{q^j}{1-xq^j}+\frac{q^j/x^2}{1-q^j/x})) = 0$. Finally,

$\lim_{z\to 0}\frac{d}{dz}\ln[(x;q)(q/x;q)/z] = \lim_{z\to 0}\frac{dx}{dz}\cdot\lim_{x\to 1}\frac{d}{dx}\ln[(x;q)(q/x;q)/z] =$
$2\pi i\cdot\lim_{x\to 1}\frac{d}{dx}\ln[(x;q)(q/x;q)/z] = \pi i$, so that $C_1(\tau) = -i\pi$. Using the well-known value

$\theta_1'(0,q) = 2q^{1/8}\prod_{j=1}^{\infty}(1-q^j)^3$, see (21) but remind that now $q = e^{2\pi i\tau}$, we find

$C(\tau) = ie^{\pi i\tau/4}(q;q)$ and thus restore (cf. e.g. [21])

$$\theta_1(\pi z \mid \tau) = i e^{\pi i(\tau/4 - z)}(x;q)(q/x;q)(q;q) \qquad (47),$$

where $(q;q) = \prod_{j=1}^{\infty}(1-q^j)$. Similar relations can be easily established for other theta-functions, including the best known [12]:

$$\theta_3(\pi z, q) = \prod_{n=1}^{\infty}(1-q^n)(1+q^{n-1/2}x)(1+q^{n-1/2}/x) \qquad (48);$$

remind that now $q = e^{2\pi i \tau}$ and $\theta_3(0,q) = \prod_{j=1}^{\infty}(1-q^j)(1+q^{j-1/2})$ [12].

## 5. Conclusions

We showed how the generalized Littlewood theorem concerning contour integrals of the logarithm of analytical function can be used in the inverse sense readily enabling to establish identities between different functions, and their transformation rules. As of now, the most interesting applications are realized for elliptical functions. Certainly, they are mostly known, although here we obtain these identities in simple and "regular" fashion - and this leads to the transparent and "standardized" forms; see, for example our Remark 3 concerning the Landen transformations, and Remark 4.

We sincerely hope that numerous other applications of this approach will be found. In saying so, we mean not only "almost evident" further applications for such generalizations of the function considered as, for example, $q$-gamma function (see e.g. [22, 23]) or elliptic gamma – function (see e.g. [21]), but also something which today is difficult to anticipate.

**Finding:** This research received no external funding.

**Data Availability Statement:** Not applicable.

**Conflicts of Interest:** The author declares no conflict of interest.